\documentclass[12pt]{article} 
\usepackage{myart2k}

\title{On the Extension of B. Sz.-Nagy's Dilation Theorem \\
to Linear Pencils of Operators}   

\author{Dmitriy S. Kalyuzhniy}

\date{}

\newcommand{\nspace}[2]{\ensuremath{{\mathbb{#1}}^{#2}}}
\newcommand{\Hspace}[1]{\ensuremath{\mathfrak{#1}}}

\def\ifundefined#1{\expandafter\ifx\csname#1\endcsname\relax}
\providecommand{\comment}[1]{}

\usepackage{amsfonts}
\providecommand{\tthdump}[1]{#1}
\tthdump{}
\newcommand{\Cliff}[2][\comment]{{\ensuremath{%
\mathcal{C}\kern-0.18em\ell(#1,#2)}}}

\ifundefined{qed}
    \DeclareMathSymbol{\qed}{0}{AMSa}{"03}
\fi
\providecommand{\eqref}[1]{\textup{(\ref{#1})}}

\providecommand{\href}[2]{#2}

\begin{document}
\maketitle
\vspace{-1cm}
\begin{abstract}
\noindent
The explicit constructions of minimal isometric, and minimal unitary
dilations of an arbitrary linear pencil of operators $T(\lambda
)=T_0+\lambda T_1$ consisting of contractions on a separable Hilbert
space for $|\lambda |=1$, which generalize the classical constructions
(the case $T_1=0$), are presented. In contrast to the classical
case these dilations are essentially non-unique.
\end{abstract}

\section{Introduction}
The classical Sz.-Nagy dilation theorem \cite{Sz.-N1} asserts that any
contractive linear operator $T$ on a Hilbert space $\Hspace{H}$ has a
unitary dilation, i.e. a unitary operator $U$ on some Hilbert space
$\Hspace{K}\supset\Hspace{H}$ such that
\begin{displaymath}
\forall n\in\mathbb{Z}_+ \quad T^n=P_\Hspace{H}U^n|\Hspace{H}
\end{displaymath}
(here $P_\Hspace{H}$ denotes the orthogonal projector onto
$\Hspace{H},\ A|\Hspace{H}$ denotes the restriction of an operator $A$
onto $\Hspace{H}$); moreover, this unitary dilation $U$ can be chosen
minimal
(in the sense of natural partial order in the set of all unitary
dilations of $T$), that is equivalent to the following:
\begin{displaymath}
\Hspace{K}=\bigvee_{n=-\infty }^{\infty }U^n\Hspace{H}
\end{displaymath}
(here $\bigvee_{n}\Hspace{L}_n$ denotes the closure of the linear span
of subsets $\Hspace{L}_n$ in $\Hspace{K}$); the minimal unitary dilation
$U$ of a contraction $T$ is unique up to unitary equivalence.

There is a quantity of generalizations of this theorem to commutative
families of contractions (see \cite{Sz.-NF}, \cite{Pa},
\cite{Ath}, \cite{Arv3} for the bibliography), and noncommutative
families of contractions (e.g.
\cite{Arv2}, \cite{Bu}, \cite{Fr}, \cite{Po2}). In the present paper we
obtain the
extension of the Sz.-Nagy dilation theorem (in the existence part) to
linear
pencils of operators.

A linear pencil $\sum_{k=1}^Nz_k\widetilde{T}_k$ of bounded linear
operators on a Hilbert space $\widetilde{\Hspace{H}}$ is called a
\emph{dilation of a linear pencil $\sum_{k=1}^Nz_kT_k$ of bounded linear
operators} on a Hilbert space $\Hspace{H}$ if
$\widetilde{\Hspace{H}}\supset\Hspace{H}$, and each of three equivalent
conditions

\emph{(i)} $\forall z=(z_1,\ldots ,z_N)\in\nspace{C}{N}, \forall
n\in\mathbb{Z}_+ \quad
\left(\sum_{k=1}^Nz_kT_k\right)^n=P_\Hspace{H}\left(\sum_{k=1}^Nz_k
\widetilde{T}_k\right)^n|\Hspace{H}$,

\emph{(ii)} $\forall \zeta=(\zeta_1,\ldots ,\zeta_N)\in\nspace{T}{N},
\forall n\in\mathbb{Z}_+ \quad
\left(\sum_{k=1}^N\zeta_kT_k\right)^n=P_\Hspace{H}\left(\sum_{k=1}^N
\zeta_k\widetilde{T}_k\right)^n|\Hspace{H}$,

\emph{(iii)} $\forall t\in\nspace{Z}{N}_+ \quad
T^t=P_\Hspace{H}\widetilde{T}^t|\Hspace{H}$

\noindent is fulfilled; here $\nspace{T}{N}:=\{ \zeta
=(\zeta_1,\ldots,\zeta_N)\in\nspace{C}{N}:|\zeta_k|=1, k=1,\ldots ,N\}$
is the unit $N$-fold torus, $\nspace{Z}{N}_+:=\{ t=(t_1,\ldots
,t_N)\in\nspace{Z}{N}:t_k\geq 0, k=1,\ldots ,N\}$ is the discrete
positive octant, $\forall t\in\nspace{Z}{N}_+\quad T^t$ is the
\emph{$t$-th symmetrized multipower of the $N$-tuple $T=(T_1,\ldots
,T_N)$ of operators}, e.g. for $t=(1,2,0,\ldots ,0)\quad
T^t=(T_1T_2^2+T_2T_1T_2+T_2^2T_1)/3$ (obviously, for a commutative
$N$-tuple $T$ this is a usual multipower: $T^t=\prod_{k=1}^NT_k^{t_k}$).
If,
moreover, for  any $\{ j_k\}_1^n\subset\{ 1,\ldots ,N\}
$\begin{equation}     \label{eq:uniform-dil}
T_{j_1}\cdots
T_{j_n}=P_\Hspace{H}\widetilde{T}_{j_1}\cdots\widetilde{T}_{j_n}|
\Hspace{H}
\end{equation}
holds then a pencil $\sum_{k=1}^Nz_k\widetilde{T}_k$ is said  to be a
\emph{uniform dilation of a pencil $\sum_{k=1}^Nz_kT_k$}. If for each
$\zeta\in\nspace{T}{N}$ the operator $\sum_{k=1}^N\zeta_kT_k$ is
contractive (resp., isometric, unitary) then we shall
refer to the set of operators $\sum_{k=1}^N\zeta_kT_k,\
\zeta\in\nspace{T}{N}$, as a \emph{contractive} (resp., \emph{isometric,
unitary}) \emph{linear pencil}. In case when a pencil
$\sum_{k=1}^N\zeta_kT_k$ is contractive, and its dilation
$\sum_{k=1}^N\zeta_k\widetilde{T}_k$ is an isometric (resp., unitary)
pencil, the latter is said to be an \emph{isometric} (resp.,
\emph{unitary}) \emph{dilation} of $\sum_{k=1}^N\zeta_kT_k$. Contractive
linear pencils appear as pencils of main operators of multiparametric
dissipative linear stationary dynamical scattering systems (see
\cite{K2}, \cite{K3}). It was proved in \cite{K3} that a contractive
linear pencil $\sum_{k=1}^N\zeta_kT_k$ on a separable Hilbert space
allows a unitary dilation if and only if for any $N$-tuple
$C=(C_1,\ldots ,C_N)$ of commuting contractions on a common separable
Hilbert space
\begin{displaymath}
\|\sum_{k=1}^NC_k\otimes T_k\|\leq 1.
\end{displaymath}
For $N=1$ this condition is, obviously, always fulfilled. For $N=2$ it
is also always fulfilled (this follows from \cite{An}). For $N\geq 3$
this condition, in general, fails \cite{K1}. Thus, in the cases $N=1$
and $N=2$ a unitary dilation of a given contractive linear pencil is
always exists. Since in the case $N=1$ the structure of minimal unitary
dilation is well known (the Sz.-Nagy dilation theorem) we shall
concentrate
our attention on the case $N=2$.

It is convenient for the sequel to consider nonhomogeneous linear
pencils of operators $T(\lambda ):=T_0+\lambda T_1,\
\lambda\in\mathbb{T}$, instead of homogeneous ones $T_{\zeta }
:=\zeta_0T_0+\zeta_1T_1,\ \zeta =(\zeta_0,\zeta_1)\in\nspace{T}{2}$. It
is clear that $T(\lambda ),\ \lambda\in\mathbb{T}$, is a contractive
(resp., isometric, unitary) pencil if and only if
$T_{\zeta },\ \zeta\in\nspace{T}{2}$, is a contractive
(resp., isometric, unitary) pencil. The definition of dilation is
reformulated as follows. A linear pencil $\widetilde{T}(\lambda )$ of
operators on a Hilbert space $\widetilde{\Hspace{H}}$ is said to be a
\emph{dilation of a linear pencil $T(\lambda )$ of operators} on a
Hilbert
space $\Hspace{H}$ if $\widetilde{\Hspace{H}}\supset\Hspace{H}$, and
\begin{equation}\label{eq:dil}
\forall\lambda\in\mathbb{T},\ \forall n\in\mathbb{Z}_+\quad T(\lambda
)^n=P_\Hspace{H}\widetilde{T}(\lambda )^n|\Hspace{H}.
\end{equation}
Note that the dilation $\widetilde{T}(\lambda )$ of a pencil $T(\lambda
)$ is called \emph{uniform} if \eqref{eq:uniform-dil} holds, and this is
equivalent to the condition
\begin{equation}\label{eq:uniform}
\forall n\in\mathbb{N}, \forall\{\lambda_j\}_1^n\subset\mathbb{T}\quad
T(\lambda_1)\cdots
T(\lambda_n)=P_\Hspace{H}\widetilde{T}(\lambda_1)\cdots\widetilde{T}
(\lambda_n)|\Hspace{H}.
\end{equation}
We shall use the term ``minimal'' for \emph{minimal isometric dilations}
(resp., \emph{minimal unitary extensions, minimal unitary dilations,
minimal uniform isometric dilations, minimal uniform unitary dilations})
in the sense of natural partial order in the set of all isometric
dilations (resp., all unitary extensions, all unitary dilations, all
uniform isometric dilations, all uniform unitary dilations) of a given
contractive linear pencil $T(\lambda )$.

In Section~\ref{sec:isom-dil} we construct a minimal isometric
dilation of an arbitrary contractive linear pencil $T(\lambda )$. This
dilation is turned out to be uniform. We also give an example of
non-uniform minimal isometric dilation, and show that both a minimal
isometric dilation and a minimal uniform isometric dilation of a
contractive linear pencil are essentially
non-unique. In Section~\ref{sec:u-dil} we construct a minimal unitary
extension of an arbitrary isometric linear pencil. Together with the
construction of a minimal isometric dilation this gives us the
construction of a minimal unitary dilation of an arbitrary contractive
linear pencil. This dilation is also turned out to be uniform. We give
also an example of non-uniform minimal unitary dilation, and show that
both a minimal unitary dilation and a minimal uniform unitary dilation
of a contractive linear pencil are essentially
non-unique. The question on the description of all minimal isometric
(resp., unitary) dilations of a contractive linear pencil is still open.

\section{Minimal isometric dilations of contractive linear
pencils}\label{sec:isom-dil}
Let $\Hspace{X}$ be a separable Hilbert space, $\mathbb{D}$ denote the
unit disk, $H_{\Hspace{X}}^2(\mathbb{D})$ denote the Hardy space
of holomorphic $\Hspace{X}$-valued functions $x$ on $\mathbb{D}$ such
that
\begin{displaymath}
\| x\|^2=(2\pi )^{-1}\sup_{0<r<1}\int_0^{2\pi }\|
x(re^{is})\|^2\,ds<\infty ,
\end{displaymath}
$[\Hspace{X},\Hspace{X}_*]$ denote the Banach space of all bounded
linear
operators from a separable Hilbert space $\Hspace{X}$ into
a separable Hilbert space $\Hspace{X}_*$. Recall (see \cite{Sz.-NF})
that a contractive holomorphic function $\theta :\mathbb{D}\to
[\Hspace{X},\Hspace{X}_*]$ is called \emph{outer} if $\overline{\theta
H^2_\Hspace{X}(\mathbb{D})}=H^2_{\Hspace{X}_*}(\mathbb{D})$
(the closure is taken in the norm of the Hilbert space  
$H^2_{\Hspace{X}_*}(\mathbb{D})$).

Let $T(\lambda )=T_0+\lambda T_1$ be a linear pencil of contractions on
a separable Hilbert space $\Hspace{H}$, i.e.
\begin{equation}\label{eq:contr}
\forall\lambda\in\mathbb{T}\quad T(\lambda )^*T(\lambda )\leq
I_\Hspace{H}
\end{equation}
where $I_\Hspace{H}$ is the identity operator on $\Hspace{H}$. Then by
the operator Fej\'{e}r--Riesz theorem (see \cite{RR}) there exist a
separable
Hilbert space $\Hspace{Y}$ and a linear outer
$[\Hspace{H},\Hspace{Y}]$-valued function $F(z)=F_0+zF_1$ such that for
boundary values $F(\lambda )=F_0+\lambda F_1$ we have:
\begin{equation}\label{eq:factor}
\forall\lambda\in\mathbb{T}\quad F(\lambda )^*F(\lambda )=I_\Hspace{H}-
T(\lambda )^*T(\lambda ).
\end{equation}
This function $F(z)$ is determined by pencil $T(\lambda )$ uniquely, up
to unitary operator factor from the left.
 Set
\begin{equation}\label{eq:k_+-space}
\Hspace{K}_+:=\left(\bigoplus_{-\infty
}^{-1}\Hspace{Y}\right)\oplus\Hspace{H},
\end{equation}
and define the operators
\begin{equation}\label{eq:canon-min-isom}
V(\lambda ):=\left[\begin{array}{cccc}
\ddots &                  &                  &             \\
       & I_{\Hspace{Y}} &                  &             \\
       &                  & I_{\Hspace{Y}} &             \\
       &                  &                  & F(\lambda ) \\
       &                  &                  & T(\lambda )
\end{array}\right]:\Hspace{K}_+\to\Hspace{K}_+\quad
(\lambda\in\mathbb{T})
\end{equation}
(here and in the sequel empty places of matrices mean zeros). It follows
from \eqref{eq:factor} that $V(\lambda )$ is an isometric linear pencil.
Let us show that $V(\lambda )$ is a uniform dilation of $T(\lambda )$.
Indeed, for any $\{\lambda_j\}_1^n\subset\mathbb{T}$ and
$h\in\Hspace{H}$ (we identify such a vector $h$ with $\mbox{col}(\ldots
,0,0,h)\in\Hspace{K}_+$) we have
\begin{eqnarray}
\lefteqn{V(\lambda_1)\cdots V(\lambda_n)h=}  \label{eq:act-v} \\
 &  \mbox{col}(\ldots
,0,0,F(\lambda_n)h,F(\lambda_{n-1})T(\lambda_n)h,\ldots
,F(\lambda_1)T(\lambda_2)\cdots T(\lambda_n)h,T(\lambda_1)\cdots
T(\lambda_n)h), \nonumber
\end{eqnarray}
and therefore we obtain
\begin{displaymath}
P_\Hspace{H}V(\lambda_1)\cdots
V(\lambda_n)|\Hspace{H}=T(\lambda_1)\cdots T(\lambda_n),
\end{displaymath}
that agrees with \eqref{eq:uniform} for $\widetilde{T}(\lambda
)=V(\lambda ),\ \lambda\in\mathbb{T}$. Thus, $V(\lambda )$ is a uniform
isometric dilation of $T(\lambda )$.
\begin{prop}\label{prop:min-isom}
The isometric dilation $\widetilde{T}(\lambda )\in
[\widetilde{\Hspace{H}}]:=[\widetilde{\Hspace{H}},\widetilde{\Hspace{H}}
],
\ \lambda\in\mathbb{T}$, of the linear pencil of contractions
$T(\lambda )\in [\Hspace{H}],\ \lambda\in\mathbb{T}$, is a minimal
isometric
dilation of $T(\lambda )$ if and only if
\begin{equation}\label{eq:min-isom}
\widetilde{\Hspace{H}}=\bigvee_{n\in\mathbb{Z}_+,\
\{\lambda_j\}_1^n\subset\mathbb{T}}\widetilde{T}(\lambda_1)\cdots
\widetilde{T}(\lambda_n)\Hspace{H}
\end{equation}
(here for $n=0$ the corresponding term is $\Hspace{H}$). If
$\widetilde{T}(\lambda )$ is a minimal uniform isometric dilation of
$T(\lambda )$ then $\widetilde{T}(\lambda )$ is a minimal isometric
dilation of $T(\lambda )$.
\end{prop}
\begin{proof}
Let $\widetilde{T}(\lambda )$ be an isometric dilation of $T(\lambda )$,
and \eqref{eq:min-isom} hold. Suppose that $\Hspace{H}'$ is a subspace
of $\widetilde{\Hspace{H}},\ \Hspace{H}'\supset\Hspace{H}$, and
$T'(\lambda
):=P_{\Hspace{H}'}\widetilde{T}(\lambda )|\Hspace{H}',\
\lambda\in\mathbb{T}$, is an isometric dilation of $T(\lambda ),\
\lambda\in\mathbb{T}$. Since for any $h'\in\Hspace{H}'$ and
$\lambda\in\mathbb{T}$
\begin{displaymath}
\|\widetilde{T}(\lambda )h'\| =\| h'\| =\| T'(\lambda )h'\| =\|
P_{\Hspace{H}'}\widetilde{T}(\lambda )h'\| ,
\end{displaymath}
we have $\widetilde{T}(\lambda )h'\in\Hspace{H}'$, and $\Hspace{H}'$ is
invariant under $\widetilde{T}(\lambda )$. Therefore,
\begin{displaymath}
\widetilde{\Hspace{H}}=\bigvee_{n,\
\{\lambda_j\}_1^n}\widetilde{T}(\lambda_1)\cdots\widetilde{T}(\lambda_n)
\Hspace{H}\subset\bigvee_{n,\
\{\lambda_j\}_1^n}\widetilde{T}(\lambda_1)\cdots\widetilde{T}(\lambda_n)
\Hspace{H}'\subset\Hspace{H}'\subset\widetilde{\Hspace{H}},
\end{displaymath}
and $\Hspace{H}'=\widetilde{\Hspace{H}}$. Thus, $\widetilde{T}(\lambda
)$ is a minimal isometric dilation of $T(\lambda )$. For the rest of
this Proposition it is sufficient to prove that if
$\widetilde{T}(\lambda )$ is a minimal uniform isometric dilation of
$T(\lambda )$ then \eqref{eq:min-isom} is true. The right-hand side of
the equality in \eqref{eq:min-isom} (denote it by $\Hspace{H}''$) is an
invariant subspace in $\widetilde{\Hspace{H}}$ under operators
$\widetilde{T}(\lambda )$ for all $\lambda\in\mathbb{T}$. If
$\widetilde{T}(\lambda )$ is a uniform isometric dilation of $T(\lambda
)$ then
\begin{equation}\label{eq:t''}
T''(\lambda ):=\widetilde{T}(\lambda )|\Hspace{H}'',\quad
\lambda\in\mathbb{T},
\end{equation}
is also a  uniform isometric dilation of $T(\lambda )$. Indeed, for any
$n\in\mathbb{N}$ and $\{\lambda_j\}_1^n\subset\mathbb{T}$ we have
\begin{displaymath}
P_\Hspace{H}T''(\lambda_1)\cdots
T''(\lambda_n)|\Hspace{H}=P_\Hspace{H}\widetilde{T}(\lambda_1)\cdots
\widetilde{T}(\lambda_n)|\Hspace{H}=T(\lambda_1)\cdots T(\lambda_n).
\end{displaymath}
Besides, \eqref{eq:t''} implies that $\widetilde{T}(\lambda )$ is a
uniform isometric dilation of $T''(\lambda )$. If $\widetilde{T}(\lambda
)$ is a minimal uniform isometric dilation of
$T(\lambda )$ then
$\Hspace{H}''=\widetilde{\Hspace{H}},\ T''(\lambda
)=\widetilde{T}(\lambda)$ for all $\lambda\in\mathbb{T}$, and the proof
is complete.
\end{proof}
Now let us show that for $\widetilde{T}(\lambda )=V(\lambda ),\
\widetilde{\Hspace{H}}=\Hspace{K}_+$, where $V(\lambda )$ is defined by
\eqref{eq:canon-min-isom}, and $\Hspace{K}_+$ is defined by
\eqref{eq:k_+-space}, the equality in \eqref{eq:min-isom} is true. From
\eqref{eq:act-v} we get for any $n\in\mathbb{N},\
\{\lambda_j\}_1^n\subset\mathbb{T}$,  and $h\in\Hspace{H}$
\begin{equation}\label{eq:gener-k_+}
V(\lambda_1)\cdots V(\lambda_n)h-V(\lambda_1)\cdots V(\lambda_{n-
1})T(\lambda_n)h=\mbox{col}(\ldots ,0,0,F(\lambda_n)h,0,\ldots ,0),
\end{equation}
with only non-zero entry $F(\lambda_n)h$ of this column vector in the
$(-n)$-th $\Hspace{Y}$'s component of $\Hspace{K}_+=(\bigoplus_{-
\infty }^{-1}\Hspace{Y})\oplus\Hspace{H}$.
Since $F(z)$ is a linear outer function, it follows from
Proposition~V.2.4 in
\cite{Sz.-NF} that for any $\lambda\in\mathbb{T}$ the lineal $F(\lambda
)\Hspace{H}$ is dense in $\Hspace{Y}$. Hence vectors of the form
\eqref{eq:gener-k_+} together with vectors from $\Hspace{H}$ are dense
in $\Hspace{K}_+$, and the desired equality
\begin{equation}\label{eq:min-isom-space}
\Hspace{K}_+=\bigvee_{n\in\mathbb{Z}_+,\
\{\lambda_j\}_1^n\subset\mathbb{T}}V(\lambda_1)\cdots
V(\lambda_n)\Hspace{H}
\end{equation}
is true.

Summing up all that was said about $V(\lambda )$ in this Section, we
obtain the following.
\begin{thm}\label{thm:min-isom}
Formulas \eqref{eq:k_+-space}--\eqref{eq:canon-min-isom} define the
minimal
isometric dilation $V(\lambda )$ of a given contractive linear pencil
$T(\lambda )$. Moreover, $V(\lambda )$ is a minimal uniform isometric
dilation of $T(\lambda )$.
\end{thm}
\begin{rem}\label{rem:particular}
It is clear that in the particular case $T_1=0$ the described
construction of minimal isometric dilation coincides with the classical
one for a contraction $T_0$ (see Section~I.4 of \cite{Sz.-NF}). In this
case $\Hspace{Y}=\overline{D_{T_0}\Hspace{H}}$ where
$D_{T_0}:=(I_\Hspace{H}-T_0^*T_0)^{1/2},\
F(\lambda )=F_0=D_{T_0}\in [\Hspace{H}, \Hspace{Y}]$ and $V(\lambda
)=V_0$ for
all $\lambda\in\mathbb{T}$.
\end{rem}
Two dilations $T'(\lambda )\in [\Hspace{H}'],\ \lambda\in\mathbb{T}$,
and $T''(\lambda )\in [\Hspace{H}''],\ \lambda\in\mathbb{T}$, of a
linear pencil $T(\lambda )\in [\Hspace{H}],\ \lambda\in\mathbb{T}$, are
said to be \emph{unitarily equivalent} if there is a unitary operator
$W:\Hspace{H}'\to\Hspace{H}''$ such that

1) $Wh=h$ for all $h\in\Hspace{H}$;

2) $\forall\lambda\in\mathbb{T}\quad T''(\lambda )=WT'(\lambda )W^{-1}$.

\noindent The following example shows that minimal isometric dilations
of a contractive linear pencil are essentially non-unique and not
necessarily uniform.
\begin{example}\label{ex:min-isom}
Consider the trivial linear pencil $T(\lambda )=0$ in
$\Hspace{H}=\mathbb{C}$. Then (see Remark~\ref{rem:particular})
$\Hspace{Y}=\mathbb{C},\ F(\lambda )=1$ for all
$\lambda\in\mathbb{T},\
\Hspace{K}_+=(\bigoplus_{-\infty}^{-
1}\mathbb{C})\oplus\mathbb{C}=\bigoplus_{-\infty }^0\mathbb{C}$, and we
obtain the following minimal uniform isometric dilation of $T(\lambda
)=0$:
\begin{equation}\label{eq:forward-shift}
V(\lambda )=\left[\begin{array}{ccc}
\ddots &   &   \\
       & 1 &   \\
       &   & 1 \\
       &   & 0 \\
\end{array}\right]:\bigoplus_{-\infty }^0\mathbb{C}\to\bigoplus_{-\infty
}^0\mathbb{C},\quad
(\lambda\in\mathbb{T})
\end{equation}
 which coincides identically with the multiplicity one forward shift
operator $S$. However, the linear pencil $V'(\lambda ):=\lambda S,\
\lambda\in\mathbb{T}$, is also a minimal
uniform isometric dilation of $T(\lambda )=0$ which is not unitarily
equivalent to the linear pencil $V(\lambda )=S$. Now set
\begin{equation}\label{eq:non-uniform}
\widetilde{V}(\lambda ):=\left[\begin{array}{cccccc}
\ddots &   &   &              &                   &                   \\
       & 1 &   &              &                   &                   \\
       &   & 1 &              &                   &                   \\
       &   &   & 1/{\sqrt{2}} & \lambda /\sqrt{2} &                   \\
       &   &   &              &                   & \lambda /\sqrt{2} \\
       &   &   &              &                   & 1/{\sqrt{2}}      \\
       &   &   &-1/{\sqrt{2}} & \lambda /\sqrt{2} &          0
\end{array}\right]:\bigoplus_{-\infty }^0\mathbb{C}\to\bigoplus_{-\infty
}^0\mathbb{C}.\quad
(\lambda\in\mathbb{T})
\end{equation}
It is verified directly that $\widetilde{V}(\lambda )$ is an isometric
linear pencil. Let us show that $\widetilde{V}(\lambda )$ is a minimal
isometric dilation of the trivial linear pencil $T(\lambda )=0$, however
is not uniform. For any $\lambda\in\mathbb{T}$ and
$h\in\Hspace{H}=\mathbb{C}$ (identified with $\mbox{col}(\ldots
,0,0,h)\in\bigoplus_{-\infty }^0\mathbb{C}$) we have
\begin{eqnarray}
\widetilde{V}(\lambda )h &=& \mbox{col}(\ldots ,0,0,\lambda
h/{\sqrt{2}},h/{\sqrt{2}},0), \nonumber \\
\widetilde{V}(\lambda)^2h &=& \mbox{col}(\ldots ,0,0,\lambda h,0,0,0),
\nonumber \\
\ldots & \ldots & \ldots \label{eq:v-tilde} \\
\widetilde{V}(
\lambda)^nh
&=& \mbox{col}(\ldots ,0,0,\underbrace{\lambda h}_{-n-1\mbox{-th}}
,0,\ldots
,0) \nonumber \\
\ldots & \ldots & \ldots . \nonumber
\end{eqnarray}
Therefore, $P_\Hspace{H}\widetilde{V}(\lambda )^nh=0=T(\lambda )^nh$ for
any $\lambda\in\mathbb{T},\ h\in\Hspace{H},\ n\in\mathbb{N}$, i.e.
\eqref{eq:dil} holds with $\widetilde{T}(\lambda )=\widetilde{V}(\lambda
)$,
and
\begin{displaymath}
\bigoplus_{-\infty }^0\mathbb{C}
=\left(\bigoplus_{n=0}^{\infty
}\widetilde{V}(1)^n\Hspace{H}\right)\oplus
\widetilde{V}(-1)\Hspace{H}\subset\bigvee_{n\in\mathbb{Z}_+,\
\{\lambda_j\}_1^n\subset\mathbb{T}}\widetilde{V}(\lambda_1)\cdots
\widetilde{V}(\lambda_n)\Hspace{H}\subset\bigoplus_{-\infty
}^0\mathbb{C}.
\end{displaymath}
Thus, the linear pencil $\widetilde{V}(\lambda )$ is a minimal isometric
dilation of the linear pencil $T(\lambda )=0$. Since for any non-zero
$h\in\Hspace{H}=\mathbb{C}$  by virtue of
\eqref{eq:non-uniform} and \eqref{eq:v-tilde} we have
\begin{displaymath}
\widetilde{V}(-1)\widetilde{V}(1)h=\widetilde{V}(-1)\mbox{col}(\ldots
,0,0, h/{\sqrt{2}}, h/{\sqrt{2}},0)=\mbox{col}(\ldots ,0,0,-h),
\end{displaymath}
we get $P_\Hspace{H}\widetilde{V}(-1)\widetilde{V}(1)h=-h\neq 0=T(-
1)T(1)h$,
and this dilation is not uniform. It is clear that a pencil
$\widetilde{V}(\lambda )$ is not unitarily equivalent both to $V(\lambda
)$ and $V'(\lambda )$ since the uniformity property of a dilation
remains under unitary equivalence transformations.
\end{example}

\section{Minimal unitary dilations of contractive linear
pencils}\label{sec:u-dil}
First we shall construct a minimal unitary extension of an arbitrary
isometric linear pencil. A unitary linear pencil $U(\lambda
)=U_0+\lambda U_1\in [\Hspace{K}],\ \lambda\in\mathbb{T}$, is said to be
a \emph{unitary extension of an isometric linear pencil} $V(\lambda
)=V_0+\lambda V_1\in [\Hspace{K}_+],\ \lambda\in\mathbb{T}$, if
$\Hspace{K}_+$ is an invariant subspace in $\Hspace{K}$ under all
operators $U(\lambda ),\ \lambda\in\mathbb{T}$, and
\begin{equation}\label{eq:ext}
\forall\lambda\in\mathbb{T}\quad V(\lambda )=U(\lambda )|\Hspace{K}_+.
\end{equation}
It is easy to see that $U(\lambda )$ is a unitary extension of an
isometric linear pencil $V(\lambda )$ if and only if $U(\lambda )$ is a
uniform unitary dilation of this pencil.

Let $V(\lambda )\in [\Hspace{K}_+],\ \lambda\in\mathbb{T}$, be an
isometric linear pencil, i.e.
\begin{displaymath}
\forall\lambda\in\mathbb{T}\quad (V_0+\lambda V_1)^*(V_0+\lambda
V_1)=I_{\Hspace{K}_+}.
\end{displaymath}
Then $V_1^*V_0=0$, i.e.
$\overline{V_0\Hspace{K}_+}\perp\overline{V_1\Hspace{K}_+}$. Let us show
that
\begin{equation}\label{eq:range}
\Hspace{V}:=\bigvee_{\lambda\in\mathbb{T}}V(\lambda )
\Hspace{K}_+=\overline{V_0\Hspace{K}_+}\oplus\overline{V_1\Hspace{K}_+}.
\end{equation}
Indeed, the inclusion ``$\subset$'' is obvious. The inclusion
``$\supset$'' follows from the Fourier representations
\begin{displaymath}
\forall k_+\in\Hspace{K}_+\quad V_jk_+=(2\pi )^{-1}\int_0^{2\pi }e^{-
ijs}V(e^{is})k_+\,ds.\quad (j=0,1)
\end{displaymath}
Set
\begin{equation}\label{eq:def-spaces}
\Hspace{L}:=\Hspace{K}_+\ominus\Hspace{V},\quad \Hspace{K}_\lambda
:=\Hspace{V}\ominus V(\lambda )\Hspace{K}_+,\quad (\lambda\in\mathbb{T})
\end{equation}
and define the unitary linear pencil
\begin{equation}\label{eq:diag}
P(\lambda ):=P_{\overline{V_0\Hspace{K}_+}}+\lambda
P_{\overline{V_1\Hspace{K}_+}}\in [\Hspace{V}],\ \lambda\in\mathbb{T}.
\end{equation}
Then $P(\lambda )V(1)\Hspace{K}_+=V(\lambda )\Hspace{K}_+$ for any
$\lambda\in\mathbb{T}$, and since $P(\lambda )$ is unitary in
$\Hspace{V}$, we have $P(\lambda )\Hspace{K}_1=\Hspace{K}_\lambda$ for
any $\lambda\in\mathbb{T}$. Set
\begin{equation}\label{eq:u-space}
\Hspace{U}:=\Hspace{K}_1\oplus\Hspace{L},
\end{equation}
\begin{equation}\label{eq:q-pencil}
Q(\lambda ):=\left[\begin{array}{cc}
P(\lambda )|\Hspace{K}_1 & 0 \\
0                        & I_\Hspace{L}
\end{array}
\right]:\Hspace{U}=\Hspace{K}_1\oplus\Hspace{L}\to\Hspace{K}_+.\quad
(\lambda\in\mathbb{T})
\end{equation}
Then $Q(\lambda )$ is an isometric linear pencil, and equalities
\begin{equation}\label{eq:factor*}
I_{\Hspace{K}_+}-V(\lambda )V(\lambda )^*=
Q(\lambda )Q(\lambda )^*,\quad (\lambda\in\mathbb{T})
\end{equation}
\begin{equation}\label{eq:orthog}
V(\lambda )^*Q(\lambda )=0 \quad (\lambda\in\mathbb{T})
\end{equation}
hold. Indeed, for any fixed $\lambda\in\mathbb{T}$ the operator
$V(\lambda )$ is isometric, hence $V(\lambda )V(\lambda )^*=P_{V(\lambda
)\Hspace{K}_+}$, and by \eqref{eq:def-spaces} $I_{\Hspace{K}_+}-
V(\lambda
)V(\lambda )^*=P_{\Hspace{K}_\lambda\oplus\Hspace{L}}$; since $Q(\lambda
)$ is also an isometry and $Q(\lambda
)\Hspace{U}=\Hspace{K}_\lambda\oplus\Hspace{L},\ Q(\lambda )Q(\lambda
)^*=P_{\Hspace{K}_\lambda\oplus\Hspace{L}}$, and \eqref{eq:factor*}
holds. Since $Q(\lambda
)\Hspace{U}=\Hspace{K}_\lambda\oplus\Hspace{L}=\Hspace{K}_+\ominus
V(\lambda )\Hspace{K}_+$, we have $Q(\lambda )\Hspace{U}\perp V(\lambda
)\Hspace{K}_+$, and \eqref{eq:orthog} follows. Define
\begin{equation}\label{eq:k-space}
\Hspace{K}:=\Hspace{K}_+\oplus\left(\bigoplus_{n=1}^\infty\Hspace{U}
\right),
\end{equation}
\begin{equation}\label{eq:u-dil}
U(\lambda ):=\left[\begin{array}{ccccc}
V(\lambda ) & Q(\lambda ) &              &              &  \\
            &             & I_\Hspace{U} &              &  \\
            &             &              & I_\Hspace{U} &  \\
            &             &              &              & \ddots
\end{array}\right]:\Hspace{K}\to\Hspace{K}.\quad (\lambda\in\mathbb{T})
\end{equation}
Since $V(\lambda )$ and $Q(\lambda )$ are isometric linear pencils and
due to \eqref{eq:factor*} and \eqref{eq:orthog} $U(\lambda )$ is a
unitary linear pencil. By \eqref{eq:u-dil} the subspace $\Hspace{K}_+$
is invariant under $U(\lambda ),\ \lambda\in\mathbb{T}$, and
\eqref{eq:ext} holds, i.e. $U(\lambda )$ is a unitary extension of the
isometric linear pencil $V(\lambda )$.
\begin{prop}\label{prop:min-u}
The unitary dilation $\widetilde{T}(\lambda )\in
[\widetilde{\Hspace{H}}],\ \lambda\in\mathbb{T}$, of a linear pencil of
contractions $T(\lambda )\in [\Hspace{H}],\ \lambda\in\mathbb{T}$, is a
minimal unitary dilation of this pencil if and only if
\begin{equation}\label{eq:min-u}
\widetilde{\Hspace{H}}=\bigvee_{n\in\mathbb{Z}_+,\
\{\lambda_j\}_1^n\subset\mathbb{T},\ \{ k_j\}_1^n\subset\{ -
1,1\}
}\widetilde{T}(\lambda_1)^{k_1}\cdots\widetilde{T}(\lambda_n)^{k_n}
\Hspace{H}
\end{equation}
(here for $n=0$ the corresponding term is $\Hspace{H}$). If
$\widetilde{T}(\lambda )$ is a minimal uniform unitary dilation of
$T(\lambda )$ then $\widetilde{T}(\lambda )$ is its minimal unitary
dilation.
\end{prop}
\begin{proof}
Let $\widetilde{T}(\lambda )$ be a unitary dilation of a linear pencil
$T(\lambda )$, and \eqref{eq:min-u} hold. Suppose that $\Hspace{H}'$
is a subspace of $\widetilde{\Hspace{H}},\
\Hspace{H}'\supset\Hspace{H}$,
and $T'(\lambda ):=P_{\Hspace{H}'}\widetilde{T}(\lambda )|\Hspace{H}',\
\lambda\in\mathbb{T}$, is a unitary dilation of a linear pencil
$T(\lambda )$. In the same way as in Proposition~\ref{prop:min-isom} we
show that $\Hspace{H}'$ is invariant under $\widetilde{T}(\lambda )$ and
$\widetilde{T}(\lambda )^*$ for all $\lambda\in\mathbb{T}$. Therefore,
\begin{displaymath}
\widetilde{\Hspace{H}}=\bigvee_{n,\ \{\lambda_j\}_1^n,\ \{
k_j\}_1^n}\widetilde{T}(\lambda_1)^{k_1}\cdots\widetilde{T}(\lambda_n)^{
k_n}\Hspace{H}\subset\bigvee_{n,\ \{\lambda_j\}_1^n,\ \{
k_j\}_1^n}\widetilde{T}(\lambda_1)^{k_1}\cdots\widetilde{T}(\lambda_n)^{
k_n}\Hspace{H}'\subset\Hspace{H}'\subset\widetilde{\Hspace{H}},
\end{displaymath}
and $\Hspace{H}'=\widetilde{\Hspace{H}}$. Thus, $\widetilde{T}(\lambda
)$ is a minimal unitary dilation of $T(\lambda )$. For the rest of this
Proposition it is sufficient to prove that if $\widetilde{T}(\lambda )$
is a minimal uniform unitary dilation of $T(\lambda )$ then
\eqref{eq:min-u} is true. The right-hand
side in \eqref{eq:min-u} (denote it by $\Hspace{H}''$) is a reducing
subspace in $\widetilde{H}$ for $\widetilde{T}(\lambda ),\
\lambda\in\mathbb{T}$. In the same way as in
Proposition~
\ref{prop:min-isom}  we can show that $T''(\lambda
):=\widetilde{T}(\lambda
)|\Hspace{H}'',\ \lambda\in\mathbb{T}$, is a uniform unitary dilation of
$T(\lambda )$. If $\widetilde{T}(\lambda )$ is a minimal uniform unitary
dilation of $T(\lambda )$ then $\Hspace{H}''=\widetilde{\Hspace{H}},\
T''(\lambda )=\widetilde{T}(\lambda )$ for all $\lambda\in\mathbb{T}$,
and the proof is complete.
\end{proof}
Now let us show that for $T(\lambda )=V(\lambda ),\
\Hspace{H}=\Hspace{K}_+,\ \widetilde{T}(\lambda )=U(\lambda )$,
and $\widetilde{\Hspace{H}}=\Hspace{K}$, where
$U(\lambda )$ and $\Hspace{K}$ are defined by
\eqref{eq:u-dil} and \eqref{eq:k-space} respectively, \eqref{eq:min-u}
is true. Let us identify vectors $k_+\in\Hspace{K}_+$ with
$\mbox{col}(k_+,0,0,\ldots )\in\Hspace{K}$. Then for any $\{
\lambda_j\}_1^n\subset\mathbb{T}$ and $k_+\in\Hspace{K}_+$ we have
\begin{eqnarray*}
\lefteqn{U(\lambda_1)^{-1}\cdots U(\lambda_n)^{-
1}k_+=} \\
& \mbox{col}(V(\lambda_1)^*\cdots
V(\lambda_n)^*k_+,Q(\lambda_1)^*V(\lambda_2)^*\cdots
V(\lambda_n)^*k_+,\ldots ,\\
& Q(\lambda_{n-
1})^*V(\lambda_n)^*k_+,Q(\lambda_n)^*k_+,0,0,\ldots ),
\end{eqnarray*}
and
\begin{eqnarray}
\lefteqn{U(\lambda_1)^{-1}\cdots U(\lambda_n)^{-1}k_+-U(\lambda_1)^{-
1}\cdots
U(\lambda_{n-1})^{-1}V(\lambda_n)^*k_+=} \nonumber  \\
& \mbox{col}(0,\ldots
,0,\underbrace{Q(\lambda_n)^*k_+}_{n-\mbox{th}},0,0,\ldots ).
\label{eq:gener-k}
\end{eqnarray}
Since by \eqref{eq:q-pencil} for any $\lambda\in\mathbb{T}$ we have
$Q(\lambda )^*\Hspace{K}_+=\Hspace{U}$, vectors of the form
\eqref{eq:gener-k} together with vectors from $\Hspace{K}_+$
fill $\Hspace{K}$,
and the desired equality
\begin{equation}\label{eq:min-u-dil}
\Hspace{K}=\bigvee_{n\in\mathbb{Z}_+,\
\{\lambda_j\}_1^n\subset\mathbb{T},\
\{ k_j\}_1^n\subset\{ -1,1\}}U(\lambda_1)^{k_1}\cdots
U(\lambda_n)^{k_n}\Hspace{K}_+
\end{equation}
is valid. Summing up all that
was said about $U(\lambda )$ in this Section, we obtain the
following.
\begin{thm}\label{thm:min-ext}
Formulas \eqref{eq:range}--\eqref{eq:q-pencil}, \eqref{eq:k-space}
and \eqref{eq:u-dil} define the minimal unitary extension $U(\lambda)$
of
a given isometric linear pencil $V(\lambda)$.
\end{thm}
\begin{rem}\label{rem:classical}
It is clear that in the particular case $V_1=0$ the
described construction of minimal unitary extension coincides with the
classical one for an isometry $V_0$ (see Section~I.2 of \cite{Sz.-NF}).
In
this case \eqref{eq:range} turns into
$\Hspace{V}=V_0\Hspace{K}_+$, \eqref{eq:def-spaces} turns into
$\Hspace{L}=\Hspace{K}_+\ominus V_0\Hspace{K}_+,\ \Hspace{K}_\lambda =\{
0\}\ (\lambda\in\mathbb{T})$, \eqref{eq:u-space} turns into
$\Hspace{U}=\Hspace{L}$, and \eqref{eq:q-pencil} turns into $Q(\lambda
)=I_\Hspace{L}:\Hspace{L}\to\Hspace{K}_+\ (\lambda\in\mathbb{T})$.
Thus, $\Hspace{U}$ coincides with the wandering generating
subspace $\Hspace{L}$ of the forward shift part of $V_0$, and $U(\lambda
)=U_0,\ \lambda\in\mathbb{T}$,
where
\begin{equation}\label{eq:class-ext}
U_0=\left[\begin{array}{cccc}
V_0 & I_\Hspace{L} &              &   \\
    &              & I_\Hspace{L} &   \\
    &              &              & \ddots
\end{array}\right]:\Hspace{K}\to\Hspace{K}
\end{equation}
is the classical minimal unitary extension of $V_0$.
\end{rem}
The following example shows that minimal unitary extensions of an
isometric linear pencil are essentially non-unique.
\begin{example}\label{ex:min-ext}
Let $V(\lambda )=S$ be a forward shift operator in
$\Hspace{K}_+=\bigoplus_{-\infty }^0\mathbb{C}$ for all
$\lambda\in\mathbb{T}$ (see \eqref{eq:forward-shift}). Then the
construction of minimal unitary extension gives (see
Remark~\ref{rem:classical}) $U(\lambda )=U_0$, where $U_0$ is
defined in \eqref{eq:class-ext} with $V_0=S,\ \Hspace{L}=\mathbb{C},\
\Hspace{K}=\bigoplus_{-\infty }^\infty \mathbb{C}$, i.e.
\begin{equation}\label{eq:ex}
U(\lambda )=U_0=\left[\begin{array}{cccccc}
\ddots  &   &          &      &   &   \\
        & 1 &          &      &   &   \\
        &   & 1        &      &   &   \\
        &   & \fbox{$0$} & 1  &   &   \\
        &   &          &      & 1 &   \\
        &   &          &      &   & \ddots
\end{array}\right]:\bigoplus_{-\infty }^\infty\mathbb{C}\to\bigoplus_{-
\infty }^\infty\mathbb{C}\quad (\lambda\in\mathbb{T})
\end{equation}
(here and in the sequel the frame distinguishes the $(0,0)$-th element
of an infinite matrix). However,
\begin{equation}\label{eq:ex'}
U'(\lambda ):=\left[\begin{array}{cccccc}
\ddots  &   &          &           &         &   \\
        & 1 &          &           &         &   \\
        &   & 1        &           &         &   \\
        &   & \fbox{$0$} & \lambda   &         &   \\
        &   &          &           & \lambda &   \\
        &   &          &           &         & \ddots
\end{array}\right]:\bigoplus_{-\infty }^\infty\mathbb{C}\to\bigoplus_{-
\infty }^\infty\mathbb{C}\quad (\lambda\in\mathbb{T})
\end{equation}
is  also a minimal unitary extension of the isometric linear pencil
$V(\lambda )=S$, which is not unitarily equivalent to the linear pencil
$U(\lambda )$.
\end{example}
\begin{prop}\label{prop:combine}
If $V(\lambda )\in [\Hspace{K}_+],\ \lambda\in\mathbb{T}$, is a minimal
isometric dilation of a contractive linear pencil $T(\lambda )\in
[\Hspace{H}],\ \lambda\in\mathbb{T}$, and $U(\lambda )\in [\Hspace{K}],\
\lambda\in\mathbb{T}$, is a minimal unitary extension of $V(\lambda )$,
then $U(\lambda )$ is a minimal unitary dilation of a pencil $T(\lambda
)$. For that, $U(\lambda )$ is a minimal uniform unitary dilation of
$T(\lambda )$ if and only if $V(\lambda )$ is a minimal uniform
isometric dilation of $T(\lambda )$.
\end{prop}
\begin{proof}
By Proposition~\ref{prop:min-isom}
the equality in \eqref{eq:min-isom-space}
is true. Since $U(\lambda )$ is a minimal uniform unitary dilation of
$V(\lambda )$ (see our remark in the beginning of this Section), by
Proposition~\ref{prop:min-u} we have \eqref{eq:min-u-dil}. Since
\eqref{eq:ext} holds, we get
\begin{eqnarray*}
\Hspace{K} &=& \bigvee_{n,\
\{\lambda_j\}_1^n,\ \{ k_j\}_1^n} U(\lambda_1)^{k_1}\cdots
U(\lambda_n)^{k_n}\Hspace{K}_+ \\
&=& \bigvee_{n,\
\{\lambda_j\}_1^n,\ \{ k_j\}_1^n} U(\lambda_1)^{k_1}\cdots
U(\lambda_n)^{k_n}\left(\bigvee_{m,\
\{\lambda_j\}_{n+1}^{n+m}}V(\lambda_{n+1})\cdots V(\lambda_{n+m})
\Hspace{H}\right) \\
&=& \bigvee_{n,\ m,\
\{\lambda_j\}_1^{n+m},\ \{ k_j\}_1^n} U(\lambda_1)^{k_1}\cdots
U(\lambda_n)^{k_n}U(\lambda_{n+1})\cdots U(\lambda_{n+m})\Hspace{H} \\
& \subset & \bigvee_{n,\
\{\lambda_j\}_1^n,\ \{ k_j\}_1^n} U(\lambda_1)^{k_1}\cdots
U(\lambda_n)^{k_n}\Hspace{H}\subset\Hspace{K},
\end{eqnarray*}
and by Proposition~\ref{prop:min-u} $U(\lambda )$ is a minimal unitary
dilation of a pencil $T(\lambda )$ (of course, a dilation of a dilation
of $T(\lambda )$ is again a dilation of $T(\lambda )$). Since for any
$\{\lambda_j\}_1^n\subset\mathbb{T}$
\begin{displaymath}
P_\Hspace{H}U(\lambda_1)\cdots
U(\lambda_n)|\Hspace{H}=P_\Hspace{H}V(\lambda_1)\cdots
V(\lambda_n)|\Hspace{H},
\end{displaymath}
$U(\lambda )$ is a uniform dilation of $T(\lambda )$ if and only if
$V(\lambda )$ is a uniform dilation of $T(\lambda )$, and the proof is
complete.
\end{proof}
Now from Theorems~\ref{thm:min-isom} and \ref{thm:min-ext}, and
Proposition~\ref{prop:combine} we obtain the following.
\begin{thm}
Formulas \eqref{eq:k_+-space}--\eqref{eq:canon-min-isom},
\eqref{eq:range}--\eqref{eq:q-pencil}, \eqref{eq:k-space} and
\begin{eqnarray}
U(\lambda ) &:=& \left[\begin{array}{cccccccc}
\ddots  & & & & & & &   \\
        & I_{\Hspace{Y}} & & & & & & \\
        & & I_{\Hspace{Y}} & & & & & \\
        & & & F(\lambda ) & P_{\Hspace{Y}}Q(\lambda ) & & & \\
        & & & \fbox{$T(\lambda )$} & P_\Hspace{H}Q(\lambda ) & & & \\
        & & & & &                   I_\Hspace{U} & &  \\
        & & & & & &                   I_\Hspace{U} &  \\
        & & & & & & &                    \ddots
\end{array}\right] \label{eq:min-unitary-dil}\\
& : & \left(\bigoplus_{-\infty }^{-1}\Hspace{Y}\right)
\oplus\Hspace{H}\oplus\left(\bigoplus_1^\infty\Hspace{U}\right)\to\left(
\bigoplus_{-\infty }^{-
1}\Hspace{Y}\right)\oplus\Hspace{H}\oplus\left(\bigoplus_1^\infty
\Hspace{U}\right)\ (=\Hspace{K})
\quad (\lambda\in\mathbb{T}) \nonumber
\end{eqnarray}
define the minimal unitary dilation $U(\lambda )$ of a given contractive
linear pencil $T(\lambda )$. Moreover, $U(\lambda )$ is a minimal
uniform unitary dilation of $T(\lambda )$.
\end{thm}
The following example shows that minimal unitary dilations of a
contractive linear pencil are essentially non-unique and not necessarily
uniform.
\begin{example}
Combining Examples~\ref{ex:min-isom} and \ref{ex:min-ext} we see that
our construction of minimal unitary dilation for the trivial linear
pencil $T(\lambda )=0\in [\mathbb{C}],\ \lambda\in\mathbb{T}$, coincides
with the classical one (see Remarks~\ref{rem:particular} and
\ref{rem:classical}) and gives identically (i.e. for all
$\lambda\in\mathbb{T}$) the multiplicity one two-sided shift operator
$U(\lambda )=U_0$ from \eqref{eq:ex}, and another minimal unitary
dilation of $T(\lambda )=0$ is $U'(\lambda )$ from \eqref{eq:ex'}. These
two minimal unitary dilations of $T(\lambda )=0$ are uniform and not
unitarily equivalent. Applying our construction of minimal unitary
extension to the non-uniform minimal isometric dilation
$\widetilde{V}(\lambda )$ of $T(\lambda )=0$,
from \eqref{eq:non-uniform}, we obtain according to
Proposition~\ref{prop:combine} the following non-uniform minimal unitary
dilation of $T(\lambda )=0$:
\begin{eqnarray*}
\widetilde{U}(\lambda ) &:=& \left[\begin{array}{cccccccccc}
\ddots  &   & & & & & & &  \\
        & 1 & & & & & & &  \\
        & & 1 & & & & & &  \\
        & & & 1/\sqrt{2} & \lambda /\sqrt{2} & & & & & \\
        & & & & & \lambda /\sqrt{2} & \lambda /\sqrt{2} & & & \\
        & & & & & 1 /\sqrt{2}       & -1/\sqrt{2}       & & & \\
        & & & -1/\sqrt{2} & \lambda /\sqrt{2} & \fbox{$0$} & & & & \\
        & & & & & & & 1 & & \\
        & & & & & & & & 1 & \\
        & & & & & & & &   & \ddots \\
\end{array}\right] \\
&:&  \bigoplus_{-\infty
}^\infty\mathbb{C}\to\bigoplus_{-\infty }^\infty\mathbb{C}.
\quad (\lambda\in\mathbb{T})
\end{eqnarray*}
It is clear that the linear pencil $\widetilde{U}(\lambda )$ is not
unitarily equivalent both to $U(\lambda )$ and $U'(\lambda )$ since the
uniformity property of a dilation remains under unitary equivalence
transformations.
\end{example}
\begin{rem}
In our construction of a minimal uniform unitary dilation of a
contractive linear pencil $T(\lambda )$ (see \eqref{eq:min-unitary-dil})
we
obtain a linear function
\begin{displaymath}
\theta (z)=\left[\begin{array}{cc}
\theta_{11}(z) & \theta_{12}(z) \\
\theta_{21}(z) & \theta_{22}(z)
\end{array}\right]:=\left[\begin{array}{cc}
F(z) & P_\Hspace{Y}Q(z) \\
T(z) & P_\Hspace{H}Q(z)
\end{array}\right]
\end{displaymath}
taking values from
$[\Hspace{H}\oplus\Hspace{U},\Hspace{Y}\oplus\Hspace{H}]$ which are
contractions for all $z\in\mathbb{D}$, and unitary operators for all
$z\in\mathbb{T}$ (i.e., a linear \emph{biinner} function). Moreover, this
function satisfies the condition
\begin{equation}\label{eq:cond}
\overline{\theta_{11}L_\Hspace{H}^2(\mathbb{T})}=
L_\Hspace{Y}^2(\mathbb{T}), \quad
\overline{\theta_{22}^*L_\Hspace{H}^2(\mathbb{T})}=
L_\Hspace{U}^2(\mathbb{T}),
\end{equation}
where $L_\Hspace{X}^2(\mathbb{T})$ denotes the Lebesgue space of 
$\Hspace{X}$-valued
square integrable functions on $\mathbb{T}$,
$\overline{\theta_{11}L_\Hspace{H}^2(\mathbb{T})}$
(resp.
$\overline{\theta_{22}^*L_\Hspace{H}^2(\mathbb{T})}$)
is the closure of the image of $L_\Hspace{H}^2(\mathbb{T})$ under 
the operator of multiplication by the boundary function $\theta_{11}(\lambda 
)$ (resp. $\theta_{22}(\lambda )^*$). The representation of an arbitrary
holomorphic contractive operator-valued function on $\mathbb{D}$ as the block
$\theta_{21}(z)$ of a biinner function $\theta (z)$ is called
\emph{Darlington's} (a $\mathcal{D}$-\emph{representation}). Let us
remark that, in general, for construction of a minimal uniform unitary 
dilation of a contractive linear pencil $T(\lambda )$ it suffices to find
a $\mathcal{D}$-representation of a function $T(z)$ with
additional requirements of linearity and fulfillment of condition
\eqref{eq:cond} on the corresponding biinner function $\theta (z)$. One
can show that Arov's general method of construction of so-called
minimal $\mathcal{D}$-representations which satisfy \eqref{eq:cond}
(see \cite{A3}), applied to a linear operator-valued function $T(z)$
which is contractive on $\mathbb{D}$, gives a linear biinner function
$\theta(z)$. Thus, minimal uniform unitary dilations of $T(\lambda )$
which are obtained in such a way deserve a special consideration.
\end{rem}

\bibliographystyle{amsplain}
\bibliography{msys}
\vskip5mm

\noindent
Department of Higher Mathematics \\
Odessa State Academy of Civil Engineering and Architecture \\
65029, Didrihson str. 4, Odessa, Ukraine \\

\end{document}